\documentclass[12pt,a4paper]{article}
\usepackage{amsthm,amsmath,amssymb,amsfonts,color, enumerate}
\usepackage[authoryear]{natbib}
 \usepackage{tikz}
\usepackage[colorlinks, allcolors=blue]{hyperref}
\usepackage{graphicx}
\usepackage{array}
\usepackage{float}
\usepackage{titling}
\usepackage{multirow}

\numberwithin{equation}{section}

\newtheoremstyle{dotless}{}{}{\itshape}{}{\bfseries}{}{ }{}
\theoremstyle{dotless}
\newtheorem{theorem}{Theorem}[section]

\newtheorem{lemma}[theorem]{Lemma}

\theoremstyle{definition}

\newtheoremstyle{definition}{}{}{}{}{\bfseries}{}{ }{}
\theoremstyle{definition}

\newtheorem{example}[theorem]{Example}

\newcommand{\R}{\mathbb{R}}

\newcommand{\bea}{\begin{eqnarray*}}
\newcommand{\eea}{\end{eqnarray*}}
\newcommand{\be}{\begin{eqnarray}}
\newcommand{\ee}{\end{eqnarray}}

\newcommand{\ww}{\omega}

\usepackage[authoryear]{natbib}
\begin{document}

\title{A note on optimal designs for estimating the slope of a  polynomial regression}
\author{\small Holger Dette \\
\small Ruhr-Universit\"at Bochum \\
\small Fakult\"at f\"ur Mathematik \\
\small 44780 Bochum, Germany \\
\small e-mail: holger.dette@rub.de\\
\and
\small Viatcheslav B. Melas \\
\small St. Petersburg State University \\
\small Department of Mathematics \\
\small St. Petersburg ,  Russia \\
{\small email: vbmelas@yandex.ru}\\
\and
\small Petr Shpilev\\
\small St. Petersburg State University \\
\small Department of Mathematics \\
\small St. Petersburg , Russia \\
{\small email: pitshp@hotmail.com}\\
}

\maketitle

\begin{abstract}
In this note we consider the  optimal design problem  for estimating  the slope of a polynomial regression with no intercept
at a given point, say $z$. In contrast to previous work, which considers symmetric design spaces we investigate  the model
 on the interval $[0, a]$ and   characterize those values of $z$, where an  explicit solution of the optimal design 
 is possible.
\end{abstract}

AMS subject classification: 62K05

Keywords and phrases: polynomial regression, slope estimation, $c$-optimal designs

\newpage

\section{Introduction}
 \label{sec1}

Consider the common polynomial regression model of degree $n$ with no intercept
\begin{equation}
 \label{1.1} Y_i = f^{\top } (x_{i} )  \theta + \varepsilon_i = 
  (x_{i},x_{i}^{2}, \dots, x_{i}^n)^\top \theta  + \varepsilon_i , \qquad i=1,\ldots,N,
\end{equation}
where $\varepsilon_1,\dots,\varepsilon_N$ denote independent  random variables with $\mathbb{E}[\varepsilon_i]
=0; $ $ {\rm Var}(\varepsilon_i)=\sigma^2 >0$,
$\theta = (\theta_{1}, \ldots , \theta_{n})^{\top} \in \mathbb{R}^n$ is a vector of unknown parameters and the explanatory variables
$x_{1}, \ldots , x_{N}$ vary in the interval $[0,a]$ for some  $a>0$.
An approximate optimal  design [in the sense of~\cite{kiefer1974}]
minimizes  an appropriate function of the (asymptotic) covariance matrix of  the statistic  $\sqrt{N} \hat \theta $,
where the  $  \hat \theta $  denotes  the least squares estimate of the parameter $\theta$
in the regression model \eqref{1.1} [see~\cite{silvey1980}
or~\cite{pukelsheim2006}].

In a recent paper  \cite{detmelshp2020}  considered model \eqref{1.1} on the 
symmetric interval $[-1,1]$ and determined explicitly the approximate  optimal design for estimating the  derivative of the regression function 
$$
{d \over dx } f^{\top } (x )  \theta  \Big |_{x=z} =  \sum_{j=1}^{n} j \theta_{j} z^{j-1}  
$$
at the point $z$, which minimizes the variance of the best linear unbiased  estimate of $ \sum_{j=1}^{n} j \theta_{j} z^{j-1}$.
The corresponding optimality criterion is a  special case of the well known
$c$-optimality   criterion  [see, for example, \cite{elfving1952,studden1968} or  \cite{pukelsheim2006}, Chapter 2].

In practice, however, polynomial regression models with no intercept are usually used on a positive interval, where $x$ corresponds, for example, to speed, concentration
or time,  and  the response at the initial  point $x=0$ is [see, for example, \cite{huangchangwong1995,lilauzhang2005}]. Therefore the goal  of this 
 note  is to provide some optimal designs for estimating the slope of  polynomial regression model with no intercept in the case where the design space is given by 
  the interval $[0,a]$.   In  Section \ref{sec2} we introduce  the basic optimal design problem and
review a geometric characterization of $c$-optimal designs. The main result can be found in  Section \ref{sec3}
where the optimal designs for estimating the slope at the point in a polynomial regression model
with no intercept are determined explicitly and the theory is illustrated by several examples.

\section{$c$-optimal designs}
\label{sec2}
\def\theequation{2.\arabic{equation}}
\setcounter{equation}{0}

Consider the regression model \eqref{1.1} on the interval $[0,a]$. Following \cite{kiefer1974} we call a discrete probability measure
 $$
\xi =
\begin{pmatrix}
x_1 & \cdots & x_m\\
\omega_1 & \cdots & \omega_m
\end{pmatrix}
$$
with support points $x_{1}, \ldots , x_{m} \in [0,a] $ and weights $\omega_1 ,  \ldots , \omega_m $ an approximate design (on the interval $[0,a]$). If $N$ observations can be taken  this means that the quantities
$N \omega_i$  are rounded to integers, say $n_{i}$, with $\sum_{i=1}^{m}n_{i} =N$  and $n_{i}$ observations are taken at each experimental condition ${x_{i}}$
 ($i=1, \ldots , m$).
 For an approximate design $\xi$ we denote by
 $$M(\xi) =\int_{[0,a]}  f(x)f^\top (x)\xi(dx)
 $$
its information matrix in the model \eqref{1.1}, where  $f(x) = (x, \ldots , x^{n} )^{\top}$ is  the  vector regression functions. The covariance matrix of the least squares estimate for the parameter $\theta$, say   $\hat \theta$,
 can be approximated (if $N \to \infty $, $n_{i} /N \to \omega_{i}$)
by $\sigma^2 /N M^{-1} (\xi)$ and an optimal design minimizes  an appropriate real valued function of the matrix  $M^{-1} (\xi)$. In this paper we are interested in designs minimizimng  the  asymptotic
 variance of the best linear unbiased estimate $c^\top \hat \theta $ of the linear combination  $c^\top \theta $ for a given vector $c \in \mathbb{R}^{n} $. To be precise, we call a design $\xi$
  $c$-optimal in the regression model \eqref{1.1}, if it minimizes the function
\[
\Phi(\xi) =
\begin{cases}
 c^\top M^{-}(\xi)c, \text{ if there exists  a vector} v \in \mathbb{R}^{d}    \text{  such that  } c = M(\xi)v; \\
 \infty, \text{ otherwise, }
\end{cases}
\]
where $M^-(\xi)$ is a generalized inverse for the matrix $M(\xi)$.  In the first case the
 design $\xi $ is called {\it admissible  for estimating the linear combination $c^\top \theta$ }  in the regression model \eqref{1.1} and the value of the quadratic form does not depend on the choice of the generalized inverse
 [see \cite{pukelsheim2006}].
The choice  $c =  f^{\prime } (z) = (1,2z ,\ldots, n z^{n-1}) ^ \top$
for some $z $ corresponds to the  minimization of the variance of the best unbiased prediction of the  derivative of the regression function $\theta^\top f(x) $ at the
point $z$. The optimal design is called {\it optimal design for estimating derivative at the point $z$} in this case.

 A useful tool for the determination of $c$-optimal designs is a geometric characterization of the $c$-optimal design and which is called Elfving's theorem in the literature 
   [see \cite{elfving1952}]. We formulate it here in a  slightly different  form, which can be directly
 used to check optimality of a given design [see \cite{DETTE2004201} for details].

\begin{theorem} \label{Elfving}
An admissible  design $\xi^{*}$  for estimating the linear combination $c^\top \theta$  with support points $x_1,  x_2,   \ldots , x_{m}  \in {\cal X} = [0,a] $ and weights
$\omega_1 , \omega_2,  \ldots , \omega_{m} $  is $c$-optimal if  and only if there exists a
vector   $p \in \mathbb{R}^{n}$  and a constant $h$ such that the following conditions are satisfied:
\begin{itemize}
	\item[(1)]
	$|p^\top f(x)|\leq1$ for all $ x\in\mathcal{X}$;
	\item[(2)]
	 $|p^\top f(x_i)   |=  1 $ for all $ i=1,2,\ldots,m$ ;
	\item[(3)]
	$c = h\sum_{i=1}^{m} f(x_i)\omega_i p^\top f(x_i)$.
\end{itemize}
Moreover, in this case we have $c^\top M^{-}(\xi^*)c=h^2$ and  the function $p^\top f(x)$  is { called extremal polynomial}.
\end{theorem}

\section{Optimal designs for estimating the slope}
\label{sec3}

For the linear model through the origin (that is $n=1$) it is easy to see  using Elfving's theorem that the
optimal design for estimating the slope is 
unique and puts all mass at the point $a$ 
 (independently of the point $z$). However, in the case $n > 1$ the situation is  more   complicated. 
By  Theorem \ref{Elfving}  it follows that the support points of the optimal design are extremal points of a 
polynomial  of  the form $p^\top f(x) = \sum_{i =1}^{n} p_{i} x^{i}$.   In fact it is possible to identify 
a candidate for this optimal polynomial explicitly. For this purpose let
\[
T_{n}(x) =\cos (n\arccos (x) )
\]
denote  the $n$th Chebyshev polynomial of the first kind
[see~\cite{szego1975}] and consider
the polynomial
\begin{equation}
\label{2.2}
  S_{n}(x) = T_n\Bigl(\frac{x}{a}(1+\cos\frac{\pi}{2n})-\cos\frac{\pi}{2n}\Bigr).
\end{equation}
\bibliographystyle{}
 It is easy to see that  $S_{n}(x)$ has exactly $n$ extremal points $s_1<s_2<\cdots<s_{n}$ on the interval $[0, a]$, which are given  by
\begin{equation}
\label{chebypoints}
  s_{i} = a\cdot\frac{\cos\frac{(i-1) \pi}{n}+\cos\frac{\pi}{2n}}{1+\cos\frac{\pi}{2n}},   ~~i=1,2,\dots, n.
\end{equation}
For the statement of our  main result we define
 $\bar L_1, \ldots , \bar L_{n}$ as the  Lagrange basis interpolation polynomials without intercept
corresponding to  the nodes $s_1, \ldots , s_n$, that is
\begin{equation}\label{3.2a}
 \bar L_i(z)=\dfrac{z\prod_{j\neq i}(z-s_j)}{s_i\prod_{j\neq i}(s_i-s_j)},
\end{equation}
and denote by   $\bar L_{i}^{\prime}$  the derivative of $\bar L_{i}$ ($i=1, \ldots , n$).

\begin{theorem}
 \label{Theorem3.2}
Consider the polynomial regression model of degree $n > 1$ with no intercept on the interval $[0,a]$.
The optimal design $\xi^*(z)$ for estimating the slope of this model at the point $z$
is supported at the points  $s_{1}, \ldots , s_{n}$
defined in \eqref{chebypoints}  if and only if  
$$
 z \in \bigcup_{i=1}^{n}(\ww_{1,i-1}, \ww_{n,i}),
 $$
where $-\ww_{1,0}=\ww_{n,n}=\infty$ and   $\ww_{i,k}$ is 
 $k$-th root of the function 
\begin{equation}  \label{weightder}
\omega_i(z)= \dfrac{|\bar L'_i(z)|}{\sum_{j=1}^{n}|\bar L'_j(z)|}~,~ i=1,\ldots,n,
\end{equation}
$k=1,\ldots,n-1,j=1,\ldots,n.$  Moreover, in this case  the weight of  the design $\xi^*(z)$
at $s_{i}$ is given  by $\omega_i(z)$ ($i=1 \ldots , n$).
\end{theorem}
To prove this Theorem we use  the following Lemma. The proof can be found in \cite{sahm1998} or in  \cite{detmelshp2020}.

\begin{lemma} \label{sahm}
Let $P_1(x)$ and $P_2(x)$ be polynomials of degree $n$ with $n$ distinct roots
$ t_{(1,1)}<t_{(1,2)}<\ldots<t_{(1,n)} $  and $ t_{(2,1)}<t_{(2,2)}<\ldots<t_{(2,n)}$, respectively.
Assume that the roots  are interlacing in the following sense:
$$
t_{(1,1)}\leq t_{(2,1)}< t_{(1,2)}\leq t_{(2,2)}<\ldots< t_{{(1,n)}}\leq t_{{(2,n)}}~,
$$
where  at least one of the inequalities $t_{{(\ell,1)}}\leq t_{{(\ell,2)}}$  ($\ell =1, \ldots n$) is strict.
Then the roots   $v_{(1,1)} \leq v_{(1,2)}   \leq  \ldots  \leq  v_{(1,n-1)}$ and  $ v_{(2,1)} \leq  v_{(2,2)}  \leq \ldots \leq  v_{(2,n-1)}$ of the
 derivatives $P_1^\prime(x) $ and $P_2^\prime(x) $  are strictly interlacing, that is
$$
v_{(1,1)}<v_{(2,1)}<\ldots<v_{(1,n-1)}<v_{(2,n-1)}.
$$
\end{lemma}

{\bf Proof of Theorem \ref{Theorem3.2}.}

We will check the optimality of the design $\xi^*(z)$ by  an application of Theorem \ref{Elfving}.
Note, that the polynomial $S_{n}(z)$  defined in \eqref{2.2} obviously satisfies to  conditions (1) and (2) of this theorem.

It now remains to characterize those values of $z$
such that the system of equations defined by  condition (3) in Theorem \ref{Elfving}  admits a solution with nonnegative weights
$\omega_{i}$ satisfying $\sum_{i=1}^{m} \omega_{i} =1$. Note that  condition  (3)   in Theorem \ref{Elfving} can be rewritten
in the form
\begin{equation} \label{cond3}
 c =f^{\prime}(z)= (1,z,\ldots,nz^{n-1})^\top =hF\beta ,
 \end{equation}
 where
 $$
  F=((s_j)^{i})_{i,j=1}^{n} = (f(s_{1}), \ldots , f(s_{n}) ) \in \R^{n  \times n}\ \hbox{and} \ \ \beta_i = \omega_i(p^\top f(s_i))
 $$
In order to investigate  the system of equations defined by   \eqref{cond3}
note that  the identity  $ F^{-1}F =I_{n}$ (here $I_{n}$ is the identity   matrix) implies
$$ e_i^\top F^{-1}f(s_j)=\delta_{ij}   ~~~(i,j=1, \ldots , n) ,$$
where $\delta_{ij}$ is the Kroneker symbol and $e_{i}=(0,\ldots , 0,1,0, \ldots , 0)^{\top} \in \R^{n}$ the $i$th unit vector.   As these equations characterize the $i$th Lagrange basis interpolation polynomial
$\bar L_i(z) =a_{i}^{T} f(z) $
with nodes
$s_{1}, \ldots , s_{n}$ we have
  $$e_i^\top F^{-1}f(z)= \bar L_i(z),\;i=1,\ldots,n.$$
Differentiating both sides of the equation with respect to $z$ yields 
 $$e_i^\top F^{-1}f'(z)= \bar L'_i(z),\;i=1,\ldots,n,$$
or equivalently
$$
F^{-1}f'(z)=(\bar L'_1(z),\ldots,\bar L'_{n}(z))^\top .$$
Therefore we obtain for the solution of \eqref{cond3}
$$
h \beta= (\bar L'_1(z),\ldots,\bar L'_{n}(z))^\top
$$
or equivalently (since ${\beta}_i = \omega_i(p^\top f(s_i))$)
\begin{equation}\label{cond4}
h  \beta_{i} = h \omega_{i} (-1)^{n-i} = \bar L'_i(z)~,~~i=1, \ldots , n.
\end{equation}
Consequently applying Lemma \ref{sahm} to the pairs of polynomials $\bar L_i(z),\bar L_{i+1}(z)$ from \eqref{3.2a} $i=1,\dots, n-1$  we 
obtain that the roots of functions $\bar L'_i(z)$ are strictly interlacing, that is
$$
\ww_{(n,1)}<\ldots<\ww_{(1,1)}<\ww_{(n,2)}<\ldots<\ww_{(1,2)}<\ldots<\ww_{(1,n-2)}<\ww_{(n,n-1)}<\ldots<\ww_{(1,n-1)}.
$$
This immediately implies that each of the functions $\bar L'_i(z)$ has only one root in the intervals $[\ww_{(n,i)},\ww_{(1,i)}]$ and has no roots in the intervals
$A_{1}, \ldots , A_{n} $, where the set $A_{i}$ is defined by
$A_{i}=(\ww_{1,i-1}, \ww_{n,i})$, $i=1,\dots, n$.
Moreover, for $z\in A_{n}$  we have 
$$
\mbox{sign}((-1)^{n-i}\bar L'_i(z))=1 ~,~~i=1, \ldots , n 
$$
(since $\mbox{sign}(\bar L'_i(z))=(-1)^{n-i}$), and for $z \in A_{j}$ 
$$
\mbox{sign}((-1)^{n-i}\bar L'_i(z))=(-1)^{n+j}~, ~~i,j=1,\ldots,n.
$$
This implies that
$$
|\bar L'_i(z)|=(-1)^{n+j}(-1)^{n-i}\bar L'_i(z) ~~(i,j=1,\ldots,n)
$$
 for $z \in A_{j}$.

The proof is now completed observing  \eqref{cond4}, which implies that 
for $z \in A_{j},$
 the weights are given by  
 $$
  \omega_{i}(z) = {  (-1)^{n-i}\bar L'_i(z) \over  h  } =  {  (-1)^{n+j} | \bar L'_i(z) | \over  h  }   ~~(i=1,\ldots,n)
  $$ 
  with
$h= (-1)^{n+j}\sum_{i=1}^{n}|\bar L'_i(z)|$.
\hfill $\Box$

\begin{figure}[t]
 	\center{\includegraphics[width=7cm,height=6cm]{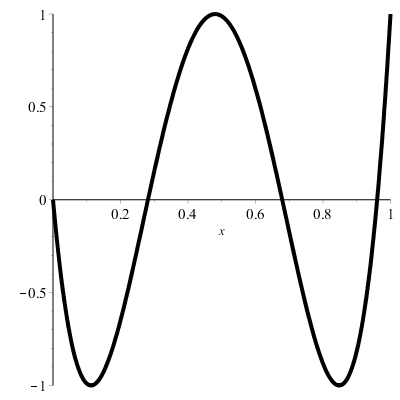}}
	\caption{\it The extremal polynomial $ S_{4} (x) $ on the interval $[0,1]$ ($a=1$). 	 \label{fig11}}
\end{figure}

\begin{figure}[t]
 	\center{\includegraphics[width=15cm]{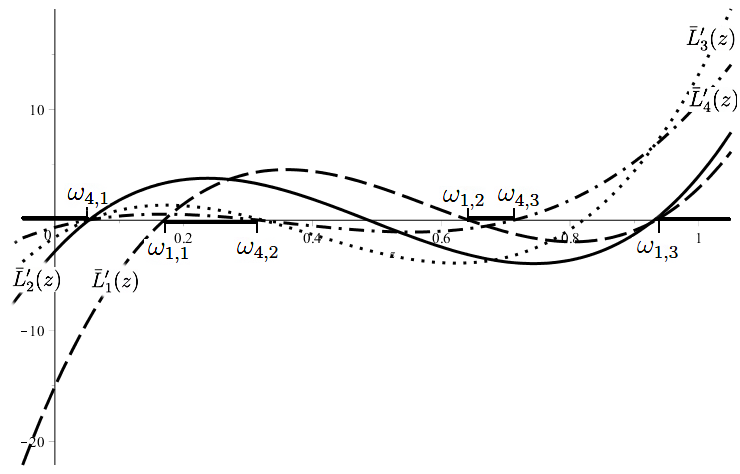}}
	\caption{\it The  functions $\bar{L}_{i}^{\prime} (z),$ $z \in \mathbb{R}$ for $i=1, 2,3 , 4$  ($n=4,\ a=1$). \label{fig1}}
\end{figure}

\bigskip

\begin{example}
\label{exam3}
{\rm In this example we illustrate potential applications of Theorem \ref{Theorem3.2} determining optimal designs for estimating the slope of a polynomial regression with no intercept
on the interval $[0,1]$.

We start with the case  of a quadratic regression model, that is $n=2$. Here  the  extremal points in  \eqref{chebypoints} are
given by  $s_{1}=\sqrt{2}-1, s_{2}=1$ and the  derivatives of the polynomials in \eqref{3.2a} are calculated as 
\begin{align*}
\bar L_{1}^{\prime}(z) &=-\frac{4+3\sqrt{2}}{2}(2z-1) ~,~ \bar L_{2}^{\prime}(z) = -\frac{2+\sqrt{2}}{2}(-2z+\sqrt{2}-1).
\end{align*}
The corresponding  roots of the functions \eqref{weightder} are obtained as
\begin{align*}
\ww_{1,1}=\frac{1}{2} ~,~ \ww_{2,1}=\frac{1}{2}(1-\sqrt{2})
\end{align*}
and the optimal design for estimating the slope of the polynomial regression without intercept  is supported at points $\sqrt{2}-1,1$
if and only if
$$
z \in \Big (-\infty,  \frac{1}{2}(1-\sqrt{2})  \Big)  \cup   \Big(\frac{1}{2}, \infty  \Big)
$$

As a second example we consider the cubic regression model with no intercept, that is $n=3$.  In this case the   extremal points of the polynomial $S_{3} (x) $
are  $s_{1}=3\sqrt{3}-5,$ $s_{2}= \sqrt{3}-1$, $s_{3}=1$ and the derivatives of the Lagrange interpolation polynomials in \eqref{3.2a} are given by
\begin{align*}
\bar L_{1}^{\prime}(z) &=35.490z^2-40.981z+8.6607 ~,~ \bar L_{2}^{\prime}(z) = -28.548z^2+22.767z-1.8680\\
\bar L_{3}^{\prime}(z) &=13.933z^2-8.6240z+.66745,
\end{align*}
The  roots of the roots of the functions \eqref{weightder}  are obtained as
\begin{align*}
\ww_{1,1}=0.2785 ~,~ \ww_{2,1}=0.0935 ~,~ \ww_{3,1}=0.090 ~,~\ww_{1,2}=0.8758 ~,~ \ww_{2,2}=0.7045 ~,~ \ww_{3,2}=0.528
\end{align*}
By Theorem \ref{Theorem3.2} 
 the optimal design for estimating the slope of the polynomial regression without intercept at the point $z$ is supported at the points   $\{3\sqrt{3}-5, \sqrt{3}-1, 1\}$
if and only if
$$
z \in (-\infty,  0.090) \cup  (0.2785, 0.528)  \cup  (0.8758, \infty)
$$
Finally we consider model \eqref{1.1} with  $n=4$, where the extremal polynomial $S_{4}(x)$ is displayed in Figure \ref{fig11}.
The corresponding extremal points 
are given by
$s_{1}=0.1127$, $s_{2}=0.4802, $, $s_{3}=0.8477$ and $s_{1}= 1$.
The  derivatives of the polynomials in \eqref{3.2a} are calculated as 
\begin{align*}
\bar L_{1}^{\prime}(z) &=-148.08z^3+258.55z^2-128.47z+15.072,\\
\bar L_{2}^{\prime}(z) &=118.63z^3-174.42z^2+62.631z-2.8327\\
\bar L_{3}^{\prime}(z) &=-114.72z^3+137.04z^2-37.110z+1.5517,\\
\bar L_{4}^{\prime}(z) &=56.968z^3-61.552z^2+15.858z-0.65327
\end{align*}
and displayed in  Figure \ref{fig1}. The roots of the functions \eqref{weightder}  are given by 
\begin{align*}
\ww_{1,1}=0.1696 ~,~ \ww_{2,1}=0.05268 ~,~ \ww_{3,1}=0.05102 ~,~\ww_{4,1}=0.05071,\\
\ww_{1,2}=0.6432 ~,~ \ww_{2,2}=0.4872 ~,~ \ww_{3,2}=0.3232 ~,~\ww_{4,2}=0.3175,\\
\ww_{1,3}=0.9332 ~,~ \ww_{2,3}=0.9305 ~,~ \ww_{3,3}=0.8205 ~,~\ww_{4,3}=0.7123.\\
\end{align*}
Therefore, by Theorem \ref{Theorem3.2}
the optimal design for estimating the slope of the polynomial regression without intercept at the point $z$ is supported at  the  points   $0.1127, 0.4802,  0.8477, 1$
if and only if
$$
z \in (-\infty,  0.05071) \cup  (0.1696, 0.3175)   \cup  (0.6432, 0.7123) \cup  (0.9332 , \infty)
$$

}

\end{example}

\bigskip
{\bf Acknowledgements}

The  work of H. Dette has been supported 
in part by the German Research Foundation, DFG (SFB 823, Teilprojekt C2, Germany's Excellence Strategy - EXC 2092 CASA - 390781972).
The work of Viatcheslav Melas and Petr Shpilev was partly supported by Russian Foundation for Basic Research (project no. 20-01-00096).

\setlength{\bibsep}{1pt}

{}

\end{document}